\documentclass[oneside]{amsart}
\usepackage{amsmath, amssymb}
\usepackage{amsthm}
\newtheorem{theorem}{Theorem}[section]

\title{A note on Haynes-Hedetniemi-Slater Conjecture}
\author{Tomoo Yokoyama}
\address{Faculty of Marine Technology, Tokyo University of Marine Science and Technology, 2-1-6, Etchujima, Koto-ku, Tokyo, 135-8533, Japan }
\date{\today}
\email{yokoyama@ms.u-tokyo.ac.jp}
 
\begin{document}

\maketitle
 
\begin{abstract}
We notice that Haynes-Hedetniemi-Slater Conjecture is true 
(i.e. $\gamma(G) \leq \frac{\delta }{3\delta -1}n$  for every graph $G$ of size $n$ with minimum degree $\delta \geq 4$, 
where $\gamma(G)$ is the domination number of $G$).
\end{abstract}

\section{Background and Remarks}

The domination number $\gamma (G)$ of a (finite, undirected and simple) graph $G = (V, E)$ 
is the minimum cardinality of
a set $D \subseteq V$ of vertices such that every vertex in $V - D$ has a neighbour in $D$.

Ore \cite{O} proved that 
$\gamma(G) \leq \frac{1}{2}n $ for every graph $G$ of size $n$ with minimum degree $\delta \geq 1$. 
Blank \cite{B} proved that 
$\gamma(G) \leq \frac{2}{5}n $ for and all but 7 exceptional graphs $G$ of size $n$  with $\delta \geq 2$. 
Blank's result was also discovered by McCuaig and Shepherd \cite{MS}. Reed \cite{R}
proved that 
$\gamma(G) \leq \frac{3}{8}n$ for every  graph $G$ of size $n$ with $\delta \geq 3$. 
All these bounds are best possible.

T. W. Haynes, S. T. Hedetniemi and P. J. Slater \cite{HHS} 
conjectured  as follows:  
\begin{quote}
For  a graph $G$ of size $n$ with minimal degree $\delta \geq 4$, 
$\gamma(G) \leq \frac{\delta}{3\delta-1}n$. 
\end{quote}
However, Caro and Roditty proved \cite{CR1} \cite{CR2} that 
for any graph $G$  of size $n$ with minimum degree 
$\delta$, 
$\gamma(G)  \leq n(1 - \delta(\frac{1}{\delta + 1})^{1 + \frac{1}{\delta }})$. 
For $\delta \geq 7$, 
it is easy to verify that
$1 - \delta(\frac{1}{\delta + 1})^{1 + \frac{1}{\delta }} < \frac{3\delta}{ \delta -1}$,  
by using calculus.
Recently, the conjecture was confirmed for $k=4$ by M. Y. Sohn and X. Yuan 
\cite{SY}
and for $k=5$ by H. M. Xing, L. Sun and X. G. Chen \cite{XSC}. 
Moreover, H.M. Xing, L. Sun, X. G. Chen prove that 
if $G$ is a Hamiltonian graph of order $n$ with $\delta \geq 6$, 
then the conjecture was confirmed for $G$. 
Therefore, the conjecture was open for graphs with $\delta = 6$. 

However, W. E. Clark, B. Shekhtman, and S. Suen, proved \cite{CSS} that 
for any graph $G$  of size $n$ with 
%minimum degree 
$\delta$, 
$\gamma(G)  \leq n(1 - \prod_{i=1}^{\delta +1} \frac{i\delta}{i \delta + 1})$. 
For $\delta = 6$, 
it is easy to verify that
$1 - \prod_{i=1}^{\delta +1} \frac{i\delta}{i \delta + 1} < \frac{3\delta}{3 \delta -1}$,  
by using calculus. 
Indeed, 
$$1 - \prod_{i=1}^{\delta +1} \frac{i\delta}{i \delta + 1} < 0.34 < \frac{6}{17} = \frac{\delta}{3 \delta -1}$$ 
and 
\begin{align*}
1 - \prod_{i=1}^{\delta +1} \frac{i\delta}{i \delta + 1} - \frac{3\delta}{3 \delta -1}  
& =\frac{11}{17}  - \prod_{i=1}^{7} \frac{6i}{6i  + 1} \\
& = - \frac{11 \cdot 7 \cdot 13 \cdot 19 \cdot 25 \cdot 31 \cdot 37 \cdot 43 - 17  \cdot 6^7 \cdot 7!}{17  \cdot 7 \cdot 13 \cdot 19 \cdot 25 \cdot 31 \cdot 37 \cdot 43}\\
& = - \frac{534014005}{36242303825} \\
& = - 0.014... \\
& < 0
\end{align*}
Therefore the conjecture is true. 

Finally we note that 
for $\delta < 6$, 
we have 
$1 - \prod_{i=1}^{\delta +1} \frac{i\delta}{i \delta + 1} > \frac{3\delta}{3 \delta -1}$.


\begin{thebibliography}{00}
\bibitem[B73]{B}
M. Blank, 
\textit{An estimate of the external stability number of a graph without suspended vertices}
Prik. Math. Programm. Vyp 10 (1973) 3-11

\bibitem[CR85]{CR1}
Y. Caro and Y. Roditty 
\textit{On the vertex-independence number and star decomposition of graphs} 
Ars Combin. 20 (1985), 167-180 

\bibitem[CR90]{CR2}
Y. Caro, Y. Roditty 
\textit{On the vertex-independence number and star decomposition of graphs} 
Internat. J. Math. Math. Sci. 13 (1990), no. 1, 205-206

\bibitem[CSS98]{CSS}
W. E. Clark, B. Shekhtman, S. Suen, 
\textit{Upper bounds of the Domination Number of a Graph} 
Congressus Numerantium, 132 (1998), pp. 99-123

\bibitem[HHS98]{HHS}
T.W. Haynes, S.T. Hedetniemi, P.J. Slater (Eds.), 
\textit{Fundamentals of domination in graphs} 
Monographs and Textbooks in Pure and Applied Mathematics, 208. Marcel Dekker, Inc., New York, 1998. 
xii+446 pp. ISBN: 0-8247-0033-3 


\bibitem[MS89]{MS}
W. McCuaig, B. Shepherd, 
\textit{Domination in graphs with minimum degree two} 
J. Graph Theory 13 (1989) 749-762.O. Ore, 

\bibitem[O62]{O}
O. Ore, 
\textit{Theory of graphs} 
Amer. Math. Soc. Colloq. Publ. 38 (1962).

\bibitem[R96]{R}
B. Reed, Paths, 
\textit{stars and the number three} 
Combin. Probab. Comput. 5 (1996) 267-276 

\bibitem[SY09]{SY}
Y. Sohn and X. Yuan 
\textit{Domination in graphs of minimum degree four} 
J. Korean Math. Soc. 46 (2009), No. 4, pp. 759-773

\bibitem[XHPA08]{XHPA}
H.M. Xing, Hattingh, J. H. Plummer, Andrew R.
\textit{On the domination number of Hamiltonian graphs with minimum degree six} 
Appl. Math. Lett. 21 (2008), no. 10, 1037-1040 

\bibitem[XSC06]{XSC}
H.M. Xing, L. Sun, X. G. Chen, 
\textit{Domination in graphs of minimum degree five} 
Graphs Combin. 22 (2006), no. 1, 127-143 


\end{thebibliography}
\end{document}